\documentclass[12pt]{article}
\def\Tilde{\char126\relax}
\usepackage{color}

\begin{document}
\begin{center}\Large{Difference Ramsey Numbers and Issai Numbers}
\end{center}
\vskip 10pt
\begin{center}Aaron Robertson\footnote{webpage: 
www.math.temple.edu/\Tilde aaron/\\ \indent This paper is part of the
author's Ph.D. thesis under the direction of Doron Zeilberger.\\
\indent This paper was supported in part by the NSF under the
PI-ship of Doron Zeilberger.}\\
Department of Mathematics, Temple University\\ Philadelphia, PA
19122\\
email:  aaron@math.temple.edu\\
\vskip 5pt
Classification:  05D10, 05D05
\end{center}
\vskip 30pt

\begin{abstract}
\noindent
We present a recursive algorithm
for finding good lower bounds for the classical Ramsey
numbers.  Using notions from this algorithm we then give
some results for generalized Schur numbers, which we
call Issai numbers.
\end{abstract}
\vskip 20pt

\begin{center}
\textbf{Introduction}
\end{center}

We present two new ideas in this paper.  The first will be 
dealing with classical Ramsey numbers.  In this part we give
a recursive algorithm for finding so-called difference Ramsey
numbers.  Using the ideas from this first part we then
define Issai numbers, a generalization of the Schur numbers.
We give some easy 
results, values, and bounds for these Issai numbers.
\vskip 10pt
Recall that $N=R(k_1,k_2,\dots,k_r)$ is the minimal integer with
the following property:
\vskip 10pt
\hskip 20pt {\bf Ramsey Property}: {\it If we $r$-color the edges of
the complete graph on $N$
\vskip 2pt
\hskip 20pt
vertices,
then there exist $j$, $1 \leq j \leq r$, such that a monochromatic $j$-colored
\vskip 2pt
\hskip 20pt
complete graph on $k_j$ vertices is a subgraph of the $r$-colored
$K_N$.}
\vskip 10pt
To find a lower bound, $L$, for one of these Ramsey numbers, it suffices to
find an
edgewise coloring of $K_L$ which avoids the Ramsey property.  To this end, 
we will restrict our seach to the subclass of {\it difference graphs}.  
After presenting some results, we will show that the Issai numbers
are a natural consequence of the difference Ramsey numbers, and a 
natural extension of the Schur numbers.
\vskip 10pt
The Difference Ramsey Numbers part of this
article is accompanied by the Maple package {\tt AUTORAMSEY}.
It has also been translated
into Fortran77 and is available as {\tt DF.f} at the author's website.
The Issai Numbers part of this article is accompanied by the
Maple package {\tt ISSAI}.
All computer packages are available for 
download at the author's website.
\newpage
\begin{center}
\textbf{Difference Ramsey Numbers}
\end{center}

Our goal here is to find good lower bounds for the classical
Ramsey numbers.  Hence, we wish to find edgewise colorings
of complete graphs which avoid the Ramsey Property.
Our approach is to construct a recursive algorithm to find
the best possible colorings among those colorings we search.
Since searching {\it all} possible colorings of a complete
graph on any nontrivial number of vertices is not feasible
by today's computing standards, we must restict the 
class of colored graphs to be searched.  The class of
graphs we will search will be the class of {\it difference
graphs.}
\vskip 10pt
\noindent
\textbf{Definition:} {\it Difference Graph}:
\vskip 2pt
\noindent  Consider the
complete graph on $n$ vertices, $K_n$.  Number the
vertices $1$ through $n$.  Let $i<j$ be two vertices of $K_n$.
Let $B_n$ be a set of arbitrary integers between $1$ and $n-1$.
Call $B_n$ the
set of blue differences on $n$ vertices.  We now color
the edges of $K_n$ as follows:  if $j-i \in B_n$ then color the edge
connecting
$i$ and $j$ blue, otherwise color the edge red.  The resulting colored
graph will be called a {\it difference graph}.   
\vskip 10pt
Given $k$ and $l$, a difference graph with the
maximal number of vertices which avoids both a blue $K_k$ and a red
$K_l$ will be called a {\it maximal difference Ramsey graph}.  Let the number
of vertices of a maximal difference Ramsey graph be $V$.  Then we
will define the {\it difference Ramsey number}, denoted $D(k,l)$, to
be $V+1$.  Further, since
the class of difference graphs is a subclass of all two-colored
complete graphs, we have that $D(k,l) \leq R(k,l)$.  Hence, by
finding the difference Ramsey numbers, we are finding lower bounds
for the classical Ramsey numbers.
\vskip 10pt
Before we present the computational aspect of these difference Ramsey
numbers, we establish an easy result:
$D(k,l) \leq D(k-1,l) + D(k,l-1)$, which is analagous to the upper bound
derived from Ramsey's proof [GRS p. 3], does {\it not} follow from Ramsey's 
proof.

\vskip 10pt

To see this consider the difference Ramsey number $D(3,3)=6$.
Let the set of red differences be $R_6=\{1,2,4\}$ (and thus 
the set of blue differences is $B_6=\{3,5\}$).  Call this
difference graph $D_6$.  In Ramsey's proof, a vertex
$v$ is isolated.   The next step is to notice that,
regardless of the choice of $v$, the number of red edges
from $v$ to $D_6 \setminus \{v\} \geq D(2,3)=3$.
Call the graph which has 
each vertex connected to the vertex $v$ by a red edge $G$.
If $v \in {1,6}$ then $G$ has $3$ vertices, otherwise it has
$4$ vertices.  Either way, the number of vertices of $G$ 
is at least $D(2,3)=3$.
 
\vskip 10pt

In order for
Ramsey's argument to work in the difference graph situation, we must show
that $G$ is isomorphic to a difference graph. 
Assume there exists an isomorphism,
$\phi: \{1,2,3,4,5,6\} \longrightarrow \{1,2,3,4,5,6\}$, such that 
the vertex set of $G$, is mapped onto $\{1,2,3\}$ or $\{1,2,3,4\}$
(depending on the number of vertices of $G$), and the
edge coloring is preserved.  Then $\phi(G)$ would be a difference graph.
Notice now that $\phi(\{v\}) \in \{4,5,6\}$.  For any choice of
$\phi(\{v\})$ we obtain the contradiction that the difference
$1$ must be both red and blue (for different edges). Hence, no such
isomorphism can exist.  Hence we cannot use the difference Ramsey number
property to conclude that the inequality holds.

However, the difference Ramsey numbers seem to be, for small values,
quite close to the Ramsey numbers.  This may just be a case of
the Law of Small Numbers, but numerical evidence from this paper
leads us to make the following
\vskip 10pt
\noindent
{\bf Conjecture 1:} $D(k,l) \leq D(k-1,l) + D(k,l-1)$
\vskip 10pt

The set of difference graphs is a superclass of the
often searched circular (or cyclic) graphs (see the survey [CG] by Chung and
Grinstead), 
which are similarly defined.  The distinction
is that, using the notation above, for a graph to be circular we
require that if $b \in B_n$, then we must have
$n-b \in B_n$.  By removing this circular condition, we remove from
the coloring the dependence on $n$ (the number of vertices), and can
thereby construct a {\it recursive} algorithm to find the set of
maximal difference Ramsey graphs: 
\vskip 10pt
The recursive step in the algorithm is described as follows.  A
difference graph on $n$ vertices consists of $B_n$, the set of
blue differences,  and $R_n$, the set of red differences.  Thus
$B_n \cup R_n = \{1,2,3,\dots,n-1\}$.
To obtain a difference graph on $n+1$ vertices, we consider
the difference $d=n$.  If $B_n \cup \{d\}$ avoids a red clique,
then we have a difference graph on $n+1$ vertices where
$B_{n+1} = B_n \cup \{d\}$ and $R_{n+1}=R_n$.  
(Note that now $B_{n+1} \cup R_{n+1} = \{1,2,3,\dots,n\}$.)
Likewise, if $R_n \cup \{d\}$ avoids a blue clique, then we
have a different difference graph on $n+1$ vertices with
$B_{n+1}=B_n$ and $R_{n+1} = R_n \cup \{d\}$.
Hence, we have a simple recursion which is not possible with
circular graphs. (By increasing the number of vertices from $n$ to
$n+1$, a circular graph goes from being circular, to being completely
noncircular (if $b \in B_n$, then 
$n-b \not \in B_n$) ).  We can now use
our recursive algorithm to find {\it automatically} (and we must note
theoretically due to time and memory constraints, but {\bf much} less
time and memory than would be required to search all graphs) all 
maximal difference Ramsey graphs for any given $k$ and $l$.
\vskip 20pt
\begin{center} {\bf About the Maple Package} {\tt AUTORAMSEY}
\end{center}

{\tt AUTORAMSEY} is a Maple package that {\it automatically} computes the 
difference graph(s) with the maximum number of vertices 
that avoids both a blue $K_k$ and a red $K_l$.  Hence, this
package {\it automatically} finds lower bounds for the Ramsey number $R(k,l)$.
In the spirit of automation, and to take another step towards AI, {\tt
AUTORAMSEY}
can create a verification Maple program tailored to
the maximal graph(s) calculated in {\tt AUTORAMSEY} (that can be run at your
leisure) and can {\it write} a \LaTeX \, paper giving the lower bound for the
Ramsey number $R(k,l)$ along with a maximal difference graph that avoids
both a blue $K_k$ and a red $K_l$.  
\vskip 10pt
The computer generated program is a straightforward program 
that can be used to (double) check that the results
obtained in {\tt AUTORAMSEY} do indeed avoid both a blue $K_k$ and a red
$K_l$.
Further, this program can be easily altered (with instructions on how
to do so) to search two-colored complete graphs for $k$-cliques and
$l$-anticliques.   
\newpage
{\tt AUTORAMSEY} has also been translated into Fortran77 
as {\tt DF.f} to speed 
up the algorithm implementation.  The code for the translated programs 
(dependent upon the clique sizes we are trying to avoid) is
 available for download at my webpage.
\vskip 20pt
\begin{center} {\bf The Algorithm}
\end{center}

Below we will give the pseudocode which finds the maximal
difference Ramsey graph(s).  Hence, it also will find the
exact value of the difference Ramsey numbers $D(k,l)$.
Because the number of difference graphs is of order $2^n$ as
compared to $2^{n^2/2}$ for all colored graphs, the algorithm
can feasibly work on larger Ramsey numbers.

Let ${\mathcal D}_n$ be the class of difference graphs on $n$ vertices.  
Let {\tt GoodSet} be the set of difference graphs that avoid both
a blue $K_k$ and a red $K_l$.  
\vskip 20pt
 $\mathtt{Let}\; m=min(k,l)$\vskip 1pt
 $\mathtt{Find}\;  {\mathcal D}_{m-1},\; \mathtt{our\ starting\ point.}$\vskip
1pt
 $\mathtt{Set}\; GoodSet = {\mathcal D}_{m-1}.$\vskip 1pt
 $\mathtt{Set}\; j = m-1 $\vskip 1pt
 $\mathtt{WHILE}\; flag \neq 0 \;\mathtt{do}$\vskip 1pt
 \hskip 20pt$\mathtt{FOR\ i\ from\ 1\ to\ \mid GoodSet \mid \; do}$\vskip 1pt
 \hskip 20pt $\mathtt{Take}\; T \in GoodSet, \;\mathtt{where}\; T
\;\mathtt{is\ of\ the\ form}\; T=[B_j,R_j] $\vskip 1pt
 \hskip 40pt $\mathtt{where}\; B_j \;\mathtt{and}\; R_j \; \mathtt{are\ the\
blue\ and\ red\ difference\ sets\ on}\;
    j \; \mathtt{vertices}$\vskip 1pt
 \hskip 20pt $\mathtt{Consider} \; S_B = [B_j \cup \{j\}, R_j] \;
\mathtt{and}\; S_R = [B_j,R_j \cup \{j\}]$\vskip 1pt
 \hskip 40pt $\mathtt{If}\; S_B \;\mathtt{avoids\ both\ a\ blue}\;
K_k\;\mathtt{and\ a\ red}\; K_l\; \mathtt{then}$\vskip 1pt
 \hskip 60pt $\mathtt{NewGoodSet := NewGoodSet \cup S_B}$\vskip 1pt
 \hskip 40pt $\mathtt{If}\; S_R \;\mathtt{avoids\ both\ a\ blue}\;
K_k\;\mathtt{and\ a\ red}\; K_l\; \mathtt{then}$\vskip 1pt
 \hskip 60pt $\mathtt{NewGoodSet := NewGoodSet \cup S_R}$\vskip 1pt
 \hskip 20pt $\mathtt{Repeat\ FOR\ loop\ with\ a\ new\ T}$\vskip 1pt
 $\mathtt{If}\; \mid NewGoodSet \mid = 0 \;\mathtt{then\ RETURN}\; GoodSet \;
   \mathtt{and\ set} \; flag=0$\vskip 1pt
 $\mathtt{Otherwise,\ set}\; GoodSet = NewGoodSet,\; NewGoodSet = \{ \},\;
\mathtt{and}\; j = j+1$\vskip 1pt 
 $\mathtt{Repeat\ WHILE\ loop}$
\vskip 20pt
For this algorithm to be efficient we must have the subroutine which
checks whether or not a monochromatic clique is avoided be very
quick.  We use the following lemma to achieve quick results in
the Fortran77 code.  (The Maple code is mainly for separatly
checking (with a different, much slower, but more
straightforward, algorithm) the Fortran77 code for small cases.)
\vskip 10pt
\noindent
\textbf{Lemma 1:}  Define the binary operation $*$ to be
$x * y = \mid x-y \mid$.  Let $D$ be a set of 
differences.
If $D$ contains a $k$-clique, then there exists
$K \subset D$, with $\mid K \mid = k-1$, such that for all $x,y \in K$, 
$x*y \in D$.
\vskip 10pt
\noindent\underline{\it Proof:} We will prove the contrapositive.
Let
$K=\{d_1,d_2,\dots,d_{k-1}\}$.  Order and rename the elements of $K$ so
that $d_1<d_2< \dots <d_{k-1}$.  Let $v_0<v_1<\dots<v_{k-1}$ be the vertices
of a $k$-clique where $d_i=v_i-v_0$. By supposition, there exists
$I<J$ such that $d_J * d_I = d_J - d_I\not \in D$.  
This is the edge connecting $v_J$ with 
$v_I$.  Since this edge is not in $D$, $D$ contains no $k$-clique.
\vskip 10pt
By using this lemma we need only check {\it pairs} of
elements in a $k$-set, rather than constructing all possible
colorings using the $k$-set.  Further, we need not worry about
the ordering of the pairs; the operation $*$ is 
commutative.
\vskip 20pt
\begin{center} {\bf Some Results}
\end{center}

It is easy to find lower bounds for $R(k,l)$, so we must show that
the algorithm gives ``good" lower bounds.  Below are two 
tables of the difference Ramsey number results obtained so far.
The first table is of the difference Ramsey number 
values.  The second table is of the number of maximal
difference Ramsey graphs.  
If we are considering the diagonal Ramsey number $R(k,k)$,
then the number of maximal difference graphs takes into account the
symmetry of colors;  i.e. we do not count a reversal of colors as a
different difference graph.
Where lower bounds are listed we have made constraints on the
size of the set {\tt GoodSet} in the algorithm due to memory
and/or (self-imposed) time restrictions. 
\vskip 20pt
\begin{center}
{\bf Difference Ramsey Numbers}
\vskip 10pt
\begin{tabular}{|lr||c|c|c|c|c|c|c|c|c|} \hline
&$\,\,\,\,l$&$3$&$4$&$5$&$6$&$7$&$8$&$9$&$10$&$11$\\      
$k$&&&&&&&&&&\\ \hline \hline
$3$&&6&9&14&17&22&27&36&39&46\\ \hline
$4$&&&18&25&34&47&$\geq 53$&$\geq 62$&&\\ \hline
$5$&&&&42&$\geq 57$&&&&&\\ \hline
\end{tabular}
\vskip 20pt
{\bf Number of Maximal Difference Ramsey Graphs}
\vskip 5pt
\begin{tabular}{|lr||c|c|c|c|c|c|c|c|c|} \hline
&$\,\,\,\,l$&$3$&$4$&$5$&$6$&$7$&$8$&$9$&$10$&$11$\\  
$k$&&&&&&&&&&\\ \hline \hline
$3$&&1&2&3&7&13&13&4&21&6\\ \hline 
$4$&&&1&6&24&21&n/a&n/a&&\\ \hline
$5$&&&&11&n/a &&&&&\\ \hline
\end{tabular}
\end{center}
\vskip 10pt

When we compare our test results to the well known maximal Ramsey graphs for
$R(3,3)$, $R(3,4)$, $R(3,5)$, $R(4,4)$ [GG], and $R(4,5)$ [MR],
we find that the program has found the critical colorings
for all of these numbers.
The classical coloring in [GRS] for $R(3,4)$ is not a difference
graph, and hence is not found by the program.  More importantly,
however, is that it does find a difference graph on $8$
vertices that avoids both a blue $K_3$ and a red $K_4$.  Hence, for
the Ramsey numbers found by Gleason and Greenwood [GG],
and for $R(4,5)$ found by McKay and Radziszowski [MR]
we have found critical Ramsey graphs which are also difference
graphs.
\vskip 20pt
The algorithm presented above can be trivially extended to
search difference graphs with more than two colors.  The progress
made so far in this direction follows.
\vskip 20pt 
\begin{center}    
{\bf Multicolored Difference Ramsey Numbers}
\end{center}

The algorithm presented here can be applied to an arbitrary number of
colors.  The recursive step in the algorithm simply becomes the
addition of the next difference to each of the three color
set $B_n$, $R_n$, and $G_n$ ($G$ for green).  Everything 
else remains the same.  Hence, the alteration of the program to any
number of colors is a simple one.  The main hurdle encountered while
searching difference graphs of more than two colors is that the 
size of the set {\tt GoodSet} in the algorithm grows {\it very}
quickly.  In fact, the system's memory while fully searching all
difference graphs was consumed within seconds for most multicolored
difference Ramsey numbers.

\begin{center}
$D(3,3,3)=15$\\
$D(3,3,4)=30$\\
$D(3,3,5)=42$\\
$D(3,3,6)\geq 60 $\\
\end{center} 
\vskip 10pt

We note here that $D(3,3,6) \geq 60$ implies that $R(3,3,6) \geq 60$, 
which is a new result.  The previous best lower bound was $54$ [SLZL].
The coloring on $59$ vertices is cyclic, hence we need only list
the differences up to $29$:
\vskip 10pt
\noindent{\tt Color 1:} 5,12,13,14,16,20,22
\vskip 5pt
\noindent{\tt Color 2:} 10,15,19,24,26,27
\vskip 5pt
\noindent {\tt Color 3:} 1,2,3,4,6,7,8,9,11,17,18,21,23,25,28,29

\vskip 20pt
\begin{center}
{\bf Future Directions}
\end{center}

Currently the algorithm which searches for the maximal difference Ramsey
graphs is a straightforward search.  In other words, if the memory 
requirement excedes the space in the computer, the algorithm will only
return a lower bound.  In the future this algorithm should be adapted
to backtrack searches or network searching.  For a backtrack search we
would note for which difference the memory barrier is reached and then
start splitting up the searches.  This would create a tree-like stucture.
We then check all leaves on this tree and choose the maximal graph.
For network searching, the same type of backtrack algorithm would be
used except that difference branches of the tree would be sent to
different computers.  This would be much quicker, but of course would
cost much more in computer facilities.
\newpage
\begin{center}
{\bf Issai Numbers}
\end{center}

Issai Schur proved in 1916 the following theorem which is considered
the first Ramsey theorem to spark activity in Ramsey Theory.
\vskip 10pt
{\it Schur's Theorem:  Given $r$, there exists an integer $N=N(r)$ such that
any
$r$-
\vskip 1pt
coloring of the integers $1$ through $N$ must admit a
monochromatic solution to
\vskip 1pt
$x+y=z$.}
\vskip 10pt
We may extend this to the following theorem:
\vskip 10pt
\noindent
{\bf Theorem 1:}  Given $r$ and $k$, there exists an integer
$N=N(r,k)$ such that any $r$-coloring of the integers $1$ through $N$
must admit a monochromatic solution to $\sum_{i=1}^{k-1} \,\, x_i =
x_k$.
\vskip 10pt
This is not a new theorem.  In fact it is a special case of Rado's Theorem 
[GRS p. 56].  We will, however, present a simple proof which relies only on
the notions already presented in this paper.
\vskip 10pt
\noindent
\underline{{\it Proof:}}  Consider the $r$-colored difference Ramsey
number
$N=D(k,k,\dots,k)$.  Then any $r$-coloring of $K_N$ must have a monochromatic
$K_k$ subgraph.  Let the vertices of this subgraph be
$\{v_0,v_1,\dots,v_{k-1}\}$,
with the differences $d_i=v_i-v_0$.  By ordering and renaming we may assume
that $d_1<d_2<\dots<d_{k-1}$.  Since $K_k$ is monochromatic, we have that
the edges $\overline{v_{i-1}v_i}$, $i=1,2,\dots,k-1$, and
$\overline{v_{k-1}v_0}$
must all be the same color.  Since the $r$-colored $K_N$ is a difference graph
we have that $(d_{i+1}-d_i)$, $i=1,2,\dots,k-1$, $d_1$, and $d_{k-1}$ must
all be assigned the same color.  Hence we have the monochromatic solution
$d_1 + \sum_{i=1}^{k-2} \,\, (d_{i+1}-d_i) =
d_{k-1}$.
\vskip 10pt

Using this theorem we will define {\it Issai numbers}.
But first, another definition is in order.
\vskip 10pt
\noindent
{\bf Definition:  }{\it Schur $k$-tuple}
\vskip 2pt
\noindent
We will call a $k$-tuple, $(x_1,x_2,\dots,x_k)$,
a {\it Schur $k$-tuple} if
$\sum_{i=1}^{k-1} \,\, x_i = x_k$.
\vskip 10pt
In the case where $k=3$, the $3$-tuple $(x,y,x+y)$ is called a Schur triple.
In Schur's theorem the only parameter is $r$, the number of colors.
Hence, a Schur number is defined to be the minimal integer $S=S(r)$ 
such that any $r$-coloring of the integers $1$ through $S$ 
must contain a monochromatic Schur triple.  
It is known that $S(2)=5$, $S(3)=14$, and $S(4)=45$. 
The Schur numbers have been generalized in [BB] and [S] in directions 
different from what will be presented here. 
We will extend the
Schur numbers in the same fashion as the Ramsey numbers
were extended from $R(k,k)$ to $R(k,l)$.
\vskip 10pt
\noindent
{\bf Definition:  }{\it Issai Number}
\vskip 2pt
\noindent
Let $S=S(k_1,k_2,\dots,k_r)$
be the minimal integer such that any $r$-coloring of the
integers from $1$ to $S$ must have a monochromatic Schur 
$k_i$-tuple, for some $i \in \{1,2,\dots,r\}$.  $S$ will be called
an {\it Issai number}.  The existence of
these Issai numbers is trivially implied by the existence
of the difference Ramsey numbers $D(k_1,k_2,\dots,k_r)$.  In fact, we
have the following result:
\vskip 10pt
\noindent
{\bf Lemma 2:}  $S(k_1,k_2,\dots,k_r) \leq D(k_1,k_2,\dots,k_r)
- 1$
\vskip 10pt
\noindent
\underline{{\it Proof:}}  By definition, there exists a minimal
integer $N=D(k_1,k_2,\dots,k_r)$ such that any $r$-coloring of $K_N$
must contain a monochromatic $K_{k_i}$, 
for some $i \in \{1,2,\dots,r\}$.  Using the same reasoning as in the proof
of Theorem 1 and the fact that the differences in the difference
graph are $1,2,\dots,N-1$, we have the stated inequality.
\vskip 10pt
Using this new definition and notation, it is already known that
$S(3,3)=5$, $S(3,3,3)=14$, and $S(3,3,3,3)=45$.  We note here that
since $D(3,3,3)=15$ we immediately have $S(3,3,3) \leq 14$, whereas
before, since $R(3,3,3)=17$, we had only that $S(3,3,3) \leq 16$.
\vskip 10pt
Attempts to find a general bound for $S(k,l)$ have been unsuccessful.
The values below lead me to make the following seemingly trivial
conjecture:
\vskip 10pt
\noindent
{\bf Conjecture 2:} $S(k-1,l) \leq S(k,l)$
\vskip 10pt
The difficulty here is that a monochromatic Schur $k$-tuple in
no way implies the existence of a monochromatic Schur $(k-1)$-tuple.
To see this, consider the following coloring of $\{1,2,\dots,9\}$.
Color $\{1,3,5,9\}$ red, and the other integers blue.
Then we have the red Schur $4$-tuple $(1,3,5,9)$. However no
red Schur triple exists in this coloring.  
\vskip 20pt
\begin{center}
{\bf Some Issai Values and Colorings}
\end{center}

We used the Maple package {\tt ISSAI}  
to calculate the
exact values as well as an
exceptional coloring given below.  {\tt ISSAI} is
written for two colors, but can easily be extended to
any number of colors.  The value $S(3,3)=5$ has been known since
before Schur proved his theorem.  
The value $S(4,4)=11$ follows from Beutelspacher and Brestovansky in [BB], 
who more generally show that $S(k,k)=k^2-k-1$.  The remaining
values are new.
\vskip 20pt
\begin{center}
{\bf Issai Numbers}
\vskip 10pt
$$
\begin{array}{|lr||c|c|c|c|c|} \hline
&\,\,\,\,l&3&4&5&6&7\\      
k&&&&&&\\ \hline \hline
3&&5&7&11&13&\geq 17\\ \hline
4&&&11&14& &\\ \hline
\end{array}
$$
\end{center}
\vskip 10pt

The exceptional colorings found by {\tt ISSAI} are as follows.
Let $S(k,l)$ denote the minimal number such that and $2$-coloring
of the integers from 1 to $S(k,l)$ must contain either a red
Schur $k$-tuple or a blue Schur $l$-tuple.  It is enough to list
only those integers colored red:
\vskip 5pt
\begin{center}
\begin{tabular}{ll}
{\tt S(3,4)>6:}
&Red:  1,6\\
{\tt S(3,5)>10:}
&Red:  1,3,8,10\\
{\tt S(4,4)>10:}
&Red:  1,2,9,10\\
{\tt S(3,6)>12:}
&Red:  1,3,10,12\\
{\tt S(4,5)>13:}
&Red:  1,2,12,13\\
{\tt S(3,7)>16:}
&Red:  1,3,5,12,14,16\\
\end{tabular}
\end{center}
\vskip 20pt
\noindent
\textbf{Acknowledgment} 
\vskip 10pt
I would like to thank my advisor, Doron Zeilberger, for
his guidence, his support, and for sharing his mathematical
philosophies. I would also like to thank Hans Johnston for
his expertise and help with my Fortran code. Further, I would
like to thank Daniel Schaal for his help with some references. 
\vskip 20pt
\noindent
\textbf{References}
\vskip 10pt
\noindent
[BB] A. Beutelspacher and W. Brestovansky,{\it Generalized
Schur Numbers}, {\tt Lecture Notes in
Mathematics} (Springer), {\bf 969}, 1982, 30-38.
\vskip 5pt
\noindent
[C] F. Chung, {\it On the Ramsey Numbers $N(3,3, \dots, 3)$}, 
Discrete Mathematics, {\bf 5}, 1973, 317-321.
\vskip 5pt
\noindent
[CG] F.R.K. Chung and C.M. Grinstead, {\it
A Survey of Bounds for Classical Ramsey Numbers},
Journal of Graph Theory, {\bf 7}, 1983, 25-37.
\vskip 5pt
\noindent
[E] G. Exoo, {\it On Two Classical Ramsey Numbers of the
Form $R(3,n)$}, SIAM Journal of Discrete Mathematics, 
{\bf 2}, 1989, 5-11.
\vskip 5pt
\noindent
[GG] A. Gleason and R. Greenwood, 
{\it Combinatorial Relations and Chromatic Graphs},
Canadian Journal of Mathematics, {\bf 7}, 1955, 1-7.
\vskip 5pt
\noindent
[GR] C. Grinstead and S. Roberts, {\it On the Ramsey Numbers
$R(3,8)$ and $R(3,9)$}, Journal of Combinatorial Theory, 
Series B, {\bf 33}, 1982, 27-51.
\vskip 5pt
\noindent 
[GRS] R. Graham, B. Rothschild, and J. Spencer,
{\bf Ramsey Theory}, John Wiley and Sons, 1980, 74-76.
\vskip 5pt
\noindent
[GY] J.E. Graver and J. Yackel, {\it Some Graph Theoretic 
Results Associated with Ramsey's Theorem}, Journal of
Combinatorial Theory, {\bf 4}, 1968, 125-175.
\vskip 5pt
\noindent
[K] J. G. Kalbfleisch, {\it Chromatic Graphs and Ramsey's
Theorem}, Ph.D. Thesis, University of Waterloo, 1966.
\vskip 5pt
\noindent
[Rad] S. Radziszowski, {\it Small Ramsey Numbers}, Electronic 
Journal of Combinatorics, Dynamic 
Survey {\bf DS1}, 1994, 28pp.
\vskip 5pt
\noindent
[RK] S. Radziszowski and D. L. Kreher, {\it On $R(3,k)$ 
Ramsey Graphs:  Theoretical and Computational Results}, 
Journal of Combinatorial Mathematics and Combinatorial 
Computing, {\bf 4}, 1988, 207-212.
\vskip 5pt
\noindent
[S] D. Schaal, {\it On Generalized Schur Numbers},
Cong. Numer., {\bf 98}, 1993, 178-187. 
\vskip 5pt
\noindent
[SLZL] Su Wenlong, Luo Haipeng, Zhang Zhengyou, and Li Guiqing, 
{\it New Lower Bounds of}     
{\it Fifteen Classical Ramsey Numbers}, to appear in Australasian
Journal of Combinatorics.
\end{document}